\documentclass[12pt]{article}%
\usepackage{amsfonts}
\usepackage{sw20bams}
\usepackage{amsmath}
\usepackage{amssymb}
\usepackage{graphicx}%
\providecommand{\U}[1]{\protect\rule{.1in}{.1in}}
\begin{document}

\title{What Do $\sin(x)$ and $\operatorname{arcsinh}(x)$ Have in Common?}
\author{Steven Finch}
\date{November 3, 2024}
\maketitle

\begin{abstract}
N. G. de Bruijn (1958) studied the asymptotic expansion of iterates of
$\sin(x)$ with $0<x\leq\pi/2$. \ Bencherif \&\ Robin (1994) generalized this
result to increasing analytic functions $f(x)$ with an attractive fixed point
at $0$ and $x>0$ suitably small. \ Mavecha \&\ Laohakosol (2013) formulated an
algorithm for explicitly deriving required parameters. \ We review their
method, testing it initally on the logistic function $\ell(x)$, a certain
radical function $r(x)$, and later on several transcendental functions.
\ Along the way, we show how $\ell(x)$ and $r(x)$ are kindred functions; the
same is also true for $\sin(x)$ and $\operatorname{arcsinh}(x)$.

\end{abstract}

\footnotetext{Copyright \copyright \ 2024 by Steven R. Finch. All rights
reserved.}Our work builds on what we began in\ \cite{F1-sinusoid, F2-sinusoid,
F3-sinusoid}. \ The focus is on nonlinear recurrences that converge to $0$.
\ While de Bruijn's book \cite{dB-sinusoid} and Bencherif \& Robin's
\cite{BR-sinusoid} paper are well-known, an article by Mavecha \&\ Laohakosol
\cite{ML-sinusoid} is completely obscure. \ A powerful algorithm developed by
these authors gives precise results that vastly exceed what was given in
\cite{F1-sinusoid, F2-sinusoid}. \ We shall demonstrate this procedure using
two algebraic functions and subsequently using logarithm/exponential
\&\ circular/hyperbolic trigonometric examples. \ The Fresnel cosine and
Lambert $W$ functions also appear. \ The algorithm does not, however, apply to
recurrences that diverge to $\infty$ and hence outcomes for some examples from
\cite{F2-sinusoid, F3-sinusoid} still await improvement.

\section{Basics}

Consider an analytic function $f(x)$ whose Taylor series at the origin starts
as
\[
x+%
{\displaystyle\sum\limits_{m=1}^{7}}
a_{m}x^{m\,\tau+1}%
\]
where $a_{1}<0$ and $\tau\geq1$ is an integer. \ We choose the upper limit $7$
for the sake of concreteness (the variable $k$ is employed in
\cite{ML-sinusoid} and leads to some confusion). \ Note that $f(0)=0$ and
zeroes of both $f$ \& $f-\operatorname{id}$ are isolated, i.e., there is an
$\varepsilon>0$ such that both $f(x)\neq0$ \& $f(x)\neq x$ for all
$0<x<\varepsilon$. \ Further, because%
\[%
\begin{array}
[c]{ccccc}%
\lim\limits_{x\rightarrow0^{+}}\dfrac{f(x)}{x}=1>0 &  & \text{\&} &  &
\lim\limits_{x\rightarrow0^{+}}\dfrac{f(x)-x}{x^{\tau+1}}=a_{1}<0,
\end{array}
\]
we see in fact that $0<f(x)<x$ for all $0<x<\varepsilon$. \ Given an initial
value $0<x_{0}<\varepsilon$, define a sequence%
\[
x_{n}=f(x_{n-1})
\]
for $n\geq1$. \ We clearly have $0<x_{n}<x_{n-1}$ and thus $\{x_{n}\}$ is
monotone decreasing \& bounded from below. \ It follows that $L=\lim
_{n\rightarrow\infty}x_{n}$ exists and $0\leq L<\varepsilon$. \ Suppose that
$L>0$. \ By definition of $x_{n}$, we must have $L=f(L)$, which contradicts
$f(L)\neq L$. \ Therefore $L=0$. \ 

Let%
\[
(-t)_{k}=(-t)(-t-1)(-t-2)\cdots(-t-k+1)
\]
for any integer $k\geq1$ and define $%
{\textstyle\sum\nolimits_{(k,m,s)}}
$ to be the sum over all nonnegative integers $n_{1},n_{2},\ldots,n_{k}$ such
that%
\[%
\begin{array}
[c]{ccccc}%
n_{1}+2n_{2}+3n_{3}+\cdots+k\,n_{k}=m &  & \text{and} &  & n_{1}+n_{2}%
+n_{3}+\cdots+n_{k}=m-s.
\end{array}
\]
This notation will be useful in the following section.

\section{Algorithm}

Given $f(x)$ and series coefficients $a_{1},a_{2},\ldots,a_{7}$ as described
in the preceding, define%
\[%
\begin{array}
[c]{ccc}%
\lambda=-\dfrac{1}{\tau\,a_{1}} &  & \text{(recalling that }\tau\text{ is the
gap between successive powers of }x\text{).}%
\end{array}
\]
Lemma A\ of \cite{ML-sinusoid} gives, for $1\leq j\leq6$,%
\[
b_{j}=\lambda^{j+1}%
{\displaystyle\sum\limits_{s=0}^{j}}
\frac{(-\tau)_{j+1-s}}{(j+1-s)!}%
{\displaystyle\sum\limits_{(7,j+1,s)}}
\dbinom{j+1-s}{n_{1},n_{2},\ldots,n_{7}}a_{1}^{n_{1}}a_{2}^{n_{2}}\cdots
a_{7}^{n_{7}}.
\]
For example,%
\[
b_{1}=\frac{(1+\tau)a_{1}^{2}-2a_{2}}{2\tau\,a_{1}^{2}}.
\]
Lemmas B \&\ C\ of \cite{ML-sinusoid} give, for $1\leq j\leq6$,%
\[
a_{0,j}=%
{\displaystyle\sum\limits_{s=0}^{j-1}}
\frac{(-1)^{j-1-s}}{j-s}%
{\displaystyle\sum\limits_{(7,j,s)}}
\dbinom{j-s}{n_{1},n_{2},\ldots,n_{7}}1^{n_{1}}b_{1}^{n_{2}}b_{2}^{n_{3}%
}\cdots b_{6}^{n_{7}}%
\]
and, for $1\leq i<j\leq6$,%
\[
a_{i,j}=%
{\displaystyle\sum\limits_{s=0}^{j-i-1}}
\frac{(-i)_{^{j-1-s}}}{(j-i-s)!}%
{\displaystyle\sum\limits_{(7,j-i,s)}}
\dbinom{j-i-s}{n_{1},n_{2},\ldots,n_{7}}1^{n_{1}}b_{1}^{n_{2}}b_{2}^{n_{3}%
}\cdots b_{6}^{n_{7}}.
\]
As examples,%
\[%
\begin{array}
[c]{ccccc}%
a_{0,1}=1, &  & a_{0,2}=b_{1}-\dfrac{1}{2}, &  & a_{0,3}=b_{2}-b_{1}+\dfrac
{1}{3},
\end{array}
\]%
\[
a_{0,4}=-\frac{b_{1}^{2}}{2}+b_{3}-b_{2}+b_{1}-\frac{1}{4};
\]%
\[%
\begin{array}
[c]{ccccc}%
a_{i,i+1}=-i, &  & a_{i,i+2}=\dfrac{i^{2}}{2}-\left(  b_{1}-\dfrac{1}%
{2}\right)  i, &  & a_{i,i+3}=-\dfrac{i^{3}}{6}+\left(  b_{1}-\dfrac{1}%
{2}\right)  i^{2}+\left(  b_{1}-b_{2}-\dfrac{1}{3}\right)  i,
\end{array}
\]%
\[
a_{i,i+4}=\frac{i^{4}}{24}-\dfrac{1}{2}\left(  b_{1}-\dfrac{1}{2}\right)
i^{3}+\dfrac{1}{2}\left(  b_{1}^{2}+2b_{2}-3b_{1}+\frac{11}{12}\right)
i^{2}-\dfrac{1}{2}\left(  2b_{3}-b_{1}^{2}-2b_{2}+2b_{1}-\frac{1}{2}\right)
i.
\]
Immediately prior to Lemma D\ of \cite{ML-sinusoid} is $c_{0}=-b_{1}$ and, for
$1\leq i\leq5$,
\[
c_{i}=\frac{1}{i}\left(  b_{i+1}+%
{\displaystyle\sum\limits_{h=0}^{i-1}}
a_{h,i+1}c_{h}\right)  .
\]
Next, define polynomials $T_{m}=T_{m}(X)$:%
\[%
\begin{array}
[c]{ccc}%
T_{1}=X, &  & T_{2}=b_{1}X-c_{1}%
\end{array}
\]
and, for $2\leq m\leq5$,%
\begin{align*}
T_{m+1}  &  =b_{1}%
{\displaystyle\sum\limits_{s=0}^{m-1}}
\frac{(-1)^{m-1-s}}{m-s}%
{\displaystyle\sum\limits_{(m,m,s)}}
\dbinom{m-s}{n_{1},n_{2},\ldots,n_{m}}T_{1}^{n_{1}}T_{2}^{n_{2}}\cdots
T_{m}^{n_{m}}-c_{m}\\
&  -%
{\displaystyle\sum\limits_{p=1}^{m-1}}
c_{p}%
{\displaystyle\sum\limits_{q=1}^{m-p}}
\frac{(-p)_{q}}{q!}%
{\displaystyle\sum\limits_{(m,m-p,m-p-q)}}
\dbinom{q}{n_{1},n_{2},\ldots,n_{m}}T_{1}^{n_{1}}T_{2}^{n_{2}}\cdots
T_{m}^{n_{m}},
\end{align*}
i.e.,%
\[
T_{3}=-\frac{b_{1}}{2}X^{2}+(b_{1}^{2}+c_{1})X-\left(  b_{1}c_{1}%
+c_{2}\right)  ,
\]%
\[
T_{4}=\frac{b_{1}}{3}X^{3}-\left(  \frac{3}{2}b_{1}^{2}+c_{1}\right)
X^{2}+(b_{1}^{3}+3b_{1}c_{1}+2c_{2})X-(b_{1}^{2}c_{1}+c_{1}^{2}+b_{1}%
c_{2}+c_{3}),
\]

\begin{align*}
T_{5}  &  =-\frac{b_{1}}{4}X^{4}+\left(  \frac{11}{6}b_{1}^{2}+c_{1}\right)
X^{3}-\left(  3b_{1}^{3}+\frac{11}{2}b_{1}c_{1}+3c_{2}\right)  X^{2}\\
&  +(b_{1}^{4}+6b_{1}^{2}c_{1}+3c_{1}^{2}+5b_{1}c_{2}+3c_{3})X-\left(
b_{1}^{3}c_{1}+\frac{5}{2}b_{1}c_{1}^{2}+b_{1}^{2}c_{2}+3c_{1}c_{2}+b_{1}%
c_{3}+c_{4}\right)  ,
\end{align*}

\begin{align*}
T_{6}  &  =\frac{b_{1}}{5}X^{5}-\left(  \frac{25}{12}b_{1}^{2}+c_{1}\right)
X^{4}+\left(  \frac{35}{6}b_{1}^{3}+\frac{25}{3}b_{1}c_{1}+4c_{2}\right)
X^{3}\\
&  -\left(  5b_{1}^{4}+\frac{35}{2}b_{1}^{2}c_{1}+6c_{1}^{2}+13b_{1}%
c_{2}+6c_{3}\right)  X^{2}\\
&  +(b_{1}^{5}+10b_{1}^{3}c_{1}+13b_{1}c_{1}^{2}+9b_{1}^{2}c_{2}+12c_{1}%
c_{2}+7b_{1}c_{3}+4c_{4})X\\
&  -\left(  b_{1}^{4}c_{1}+\frac{9}{2}b_{1}^{2}c_{1}^{2}+2c_{1}^{3}+b_{1}%
^{3}c_{2}+7b_{1}c_{1}c_{2}+2c_{2}^{2}+b_{1}^{2}c_{3}+4c_{1}c_{3}+b_{1}%
c_{4}+c_{5}\right)  .
\end{align*}
We mentioned earlier some confusion in \cite{ML-sinusoid}: their recursive
formula for $T_{m+1}$ is problematic because it seems that $T_{m+1}$ might
appear on both left \&\ right sides of the equality. \ The difficulty is
avoided by reinterpreting their $k-1$ as $m$.\ \ As stated earlier, setting
their $k$ to be $7$ everywhere else is suitable for our purposes. \ Also, set
$\widetilde{T}_{m}(X)=T_{m}(-\tau\,X)$ for convenience.

Finally, define polynomials $P_{m}=P_{m}(X)$:%
\[%
\begin{array}
[c]{ccc}%
P_{0}=1, &  & P_{1}=X,
\end{array}
\]
and, for $2\leq m\leq6$,
\[
P_{m}=%
{\displaystyle\sum\limits_{s=0}^{m-1}}
\frac{(-\frac{1}{\tau})_{m-s}}{(m-s)!}%
{\displaystyle\sum\limits_{(6,m,s)}}
\dbinom{m-s}{n_{1},n_{2},\ldots,n_{6}}\widetilde{T}_{1}^{n_{1}}\widetilde
{T}_{2}^{n_{2}}\cdots\widetilde{T}_{6}^{n_{6}},
\]
e.g.,%
\[
P_{2}=\frac{\tau}{2}\left(  1+\frac{1}{\tau}\right)  \,X^{2}+b_{1}%
X+\frac{c_{1}}{\tau},
\]%
\[
P_{3}=\frac{\tau^{2}}{6}\left(  1+\frac{1}{\tau}\right)  \left(  2+\frac
{1}{\tau}\right)  \,X^{3}+\left(  1+\frac{3\tau}{2}\right)  b_{1}%
\,X^{2}+\left(  b_{1}^{2}+\left(  2+\frac{1}{\tau}\right)  c_{1}\right)
X+\frac{b_{1}c_{1}+c_{2}}{\tau}.
\]
These polynomials satisfy a differential-difference equation
\cite{BR-sinusoid}
\[
P_{m+1}^{\prime}=b_{1}P_{m}^{\prime}+(m\,\tau+1)P_{m}%
\]
and are fundamental in establishing the asymptotics of iterates of $f(x)$.
\ Essential theory motivating such formulas may be found in \cite{dB-sinusoid,
BR-sinusoid, ML-sinusoid}. \ It suffices to note that the transformed sequence
$y_{n}=\lambda/x_{n}^{\tau}$ approaches $\infty$ with $n$ and features $b_{j}%
$, $a_{0,j}$, $a_{i,j}$ in asymptotic expansions of $y_{n+1}-y_{n}$,
$\ln(y_{n+1}/y_{n})$, $y_{n+1}^{-i}-y_{n}^{-i}$ respectively. \ Also, $c_{i}$
appears in a Laurent series for a certain expression $\Psi(x)$, whose
functional inverse $\Psi^{-1}(y)$ has a series involving $T_{j}$. \ For
reasons of space, we omit all supporting details.

\section{Logistic}

For the logistic function $\ell(x)=x(1-x)$, we have $\tau=1$,%
\[
\{a_{m}\}_{m=1}^{7}=\left\{  -1,0,0,0,0,0,0\right\}
\]
and $\lambda=1$; consequently%
\[
\{b_{j}\}_{j=1}^{6}=\left\{  1,1,1,1,1,1\right\}  ,
\]%
\[
\{a_{0j}\}_{j=1}^{6}=\left\{  1,\frac{1}{2},\frac{1}{3},\frac{1}{4},\frac
{1}{5},\frac{1}{6}\right\}  ,
\]%
\[
\{c_{i}\}_{i=1}^{5}=\left\{  \frac{1}{2},\frac{1}{3},\frac{13}{36},\frac
{113}{240},\frac{1187}{1800}\right\}  ,
\]%
\[
T_{2}=X-\dfrac{1}{2},
\]%
\[
T_{3}=-\dfrac{1}{2}X^{2}+\dfrac{3}{2}X-\dfrac{5}{6},
\]%
\[
T_{4}=\dfrac{1}{3}X^{3}-2X^{2}+\dfrac{19}{6}X-\dfrac{13}{9},
\]%
\[
T_{5}=-\dfrac{1}{4}X^{4}+\dfrac{7}{3}X^{3}-\dfrac{27}{4}X^{2}+\dfrac{15}%
{2}X-\dfrac{2009}{720},
\]%
\[
T_{6}=\dfrac{1}{5}X^{5}-\dfrac{31}{12}X^{4}+\dfrac{34}{3}X^{3}-\dfrac{87}%
{4}X^{2}+\dfrac{3359}{180}X-\dfrac{6973}{1200}%
\]

and%
\[
P_{2}=X^{2}+X+\dfrac{1}{2},
\]%
\[
P_{3}=X^{3}+\dfrac{5}{2}X^{2}+\dfrac{5}{2}X+\dfrac{5}{6},
\]%
\[
P_{4}=X^{4}+\dfrac{13}{3}X^{3}+\dfrac{15}{2}X^{2}+\dfrac{35}{6}X+\dfrac
{61}{36},
\]%
\[
P_{5}=X^{5}+\dfrac{77}{12}X^{4}+\dfrac{101}{6}X^{3}+\dfrac{265}{12}%
X^{2}+\dfrac{515}{36}X+\dfrac{2609}{720},
\]%
\[
P_{6}=X^{6}+\dfrac{87}{10}X^{5}+\dfrac{95}{3}X^{4}+61X^{3}+65X^{2}%
+\dfrac{12977}{360}X+\dfrac{29069}{3600}.
\]
Of many possible intricacies, let us focus on the derivation of $P_{m}$,
$1\leq m\leq4$, from $\widetilde{T}_{1},\widetilde{T}_{2},\widetilde{T}%
_{3},\widetilde{T}_{4}$. \ We have%
\[
\frac{(-\frac{1}{\tau})_{m-s}}{(m-s)!}=\frac{(-1)_{m-s}}{(m-s)!}%
=\frac{(-1)(-2)\cdots(-1-(m-s)+1)}{(m-s)!}=(-1)^{m-s}%
\]
since $\tau=1$. \ The nonzero contributions to the sum underlying $P_{m}$ are%
\[%
\begin{array}
[c]{ccc}%
(-1)^{1-0}\dfrac{(1-0)!}{1!}\widetilde{T}_{1}=-\widetilde{T}_{1} &  &
[m=1,s=0;n_{1}=1]
\end{array}
\]%
\[%
\begin{array}
[c]{ccc}%
(-1)^{2-0}\dfrac{(2-0)!}{2!}\widetilde{T}_{1}^{2}=\widetilde{T}_{1}^{2} &  &
[m=2,s=0;n_{1}=2]
\end{array}
\]%
\[%
\begin{array}
[c]{ccc}%
(-1)^{2-1}\dfrac{(2-1)!}{0!1!}\widetilde{T}_{2}=-\widetilde{T}_{2} &  &
[m=2,s=1;n_{1}=0,n_{2}=1]
\end{array}
\]%
\[%
\begin{array}
[c]{ccc}%
(-1)^{3-0}\dfrac{(3-0)!}{3!}\widetilde{T}_{1}^{3}=-\widetilde{T}_{1}^{3} &  &
[m=3,s=0;n_{1}=3]
\end{array}
\]%
\[%
\begin{array}
[c]{ccc}%
(-1)^{3-1}\dfrac{(3-1)!}{1!1!}\widetilde{T}_{1}\widetilde{T}_{2}%
=2\widetilde{T}_{1}\widetilde{T}_{2} &  & [m=3,s=1;n_{1}=1,n_{2}=1]
\end{array}
\]%
\[%
\begin{array}
[c]{ccc}%
(-1)^{3-2}\dfrac{(3-2)!}{0!0!1!}\widetilde{T}_{3}=-\widetilde{T}_{3} &  &
[m=3,s=2;n_{1}=n_{2}=0,n_{3}=1]
\end{array}
\]%
\[%
\begin{array}
[c]{ccc}%
(-1)^{4-0}\dfrac{(4-0)!}{4!}\widetilde{T}_{1}^{4}=\widetilde{T}_{1}^{4} &  &
[m=4,s=0;n_{1}=4]
\end{array}
\]%
\[%
\begin{array}
[c]{ccc}%
(-1)^{4-1}\dfrac{(4-1)!}{2!1!}\widetilde{T}_{1}^{2}\widetilde{T}%
_{2}=-3\widetilde{T}_{1}^{2}\widetilde{T}_{2} &  & [m=4,s=1;n_{1}=2,n_{2}=1]
\end{array}
\]%
\[%
\begin{array}
[c]{ccc}%
(-1)^{4-2}\dfrac{(4-2)!}{0!2!}\widetilde{T}_{2}^{2}=\widetilde{T}_{2}^{2} &  &
[m=4,s=2;n_{1}=0,n_{2}=2]^{\ast}%
\end{array}
\]%
\[%
\begin{array}
[c]{ccc}%
(-1)^{4-2}\dfrac{(4-2)!}{1!0!1!}\widetilde{T}_{1}\widetilde{T}_{3}%
=2\widetilde{T}_{1}\widetilde{T}_{3} &  & [m=4,s=2;n_{1}=1,n_{2}%
=0,n_{3}=1]^{\ast}%
\end{array}
\]%
\[%
\begin{array}
[c]{ccc}%
(-1)^{4-3}\dfrac{(4-3)!}{0!0!0!1!}\widetilde{T}_{4}=-\widetilde{T}_{4} &  &
[m=4,s=3;n_{1}=n_{2}=n_{3}=0,n_{4}=1]
\end{array}
\]
and the two entries with matching $(m,s)$ are starred.

At long last, we exhibit the magnificent connection between $P_{m}$ and
asymptotics of $x_{n}=\ell(x_{n-1})$:%
\begin{equation}
x_{n}\sim\left(  \frac{\lambda}{n}\right)  ^{1/\tau}\left\{  1+%
{\displaystyle\sum\limits_{m=1}^{6}}
P_{m}\left(  -\frac{1}{\tau}\left[  b_{1}\ln(n)+C\right]  \right)  \frac
{1}{n^{m}}\right\}  \tag{A}\label{A}%
\end{equation}
as $n\rightarrow\infty$, where $C=C(x_{0})$ is a constant depending on the
initial condition. \ Performing the $X$ substitution as indicated, we obtain
\[
P_{1}(X)=-T_{1}(-X)=-\ln(n)-C
\]%
\[
P_{2}(X)=T_{1}(-X)^{2}-T_{2}(-X)=\ln(n)^{2}+(2C-1)\ln(n)+\left(  C^{2}%
-C+\frac{1}{2}\right)
\]%
\begin{align*}
P_{3}(X)  &  =-T_{1}(-X)^{3}+2T_{1}(-X)T_{2}(-X)-T_{3}(-X)\\
&  =-\ln(n)^{3}-\left(  3C-\frac{5}{2}\right)  \ln(n)^{2}-\left(
3C^{2}-5C+\frac{5}{2}\right)  \ln(n)-\left(  C^{3}-\frac{5}{2}C^{2}+\frac
{5}{2}C-\frac{5}{6}\right)
\end{align*}%
\begin{align*}
P_{4}(X)  &  =T_{1}(-X)^{4}-3T_{1}(-X)^{2}T_{2}(-X)+T_{2}(-X)^{2}%
+2T_{1}(-X)T_{3}(-X)-T_{4}(-X)\\
&  =\ln(n)^{4}+\left(  4C-\frac{13}{3}\right)  \ln(n)^{3}+\left(
6C^{2}-13C+\frac{15}{2}\right)  \ln(n)^{2}\\
&  +\left(  4C^{3}-13C^{2}+15C-\frac{35}{6}\right)  \ln(n)+\left(  C^{4}%
-\frac{13}{3}C^{3}+\frac{15}{2}C^{2}-\frac{35}{6}C+\frac{61}{36}\right)
\end{align*}
which implies%
\begin{align*}
x_{n}  &  \sim\frac{1}{n}-\frac{\ln(n)}{n^{2}}-\frac{C}{n^{2}}+\frac
{\ln(n)^{2}}{n^{3}}+(2C-1)\frac{\ln(n)}{n^{3}}+\left(  C^{2}-C+\frac{1}%
{2}\right)  \frac{1}{n^{3}}\\
&  -\frac{\ln(n)^{3}}{n^{4}}-\left(  3C-\frac{5}{2}\right)  \frac{\ln(n)^{2}%
}{n^{4}}-\left(  3C^{2}-5C+\frac{5}{2}\right)  \frac{\ln(n)}{n^{4}}-\left(
C^{3}-\frac{5}{2}C^{2}+\frac{5}{2}C-\frac{5}{6}\right)  \frac{1}{n^{4}}\\
&  +\frac{\ln(n)^{4}}{n^{5}}+\left(  4C-\frac{13}{3}\right)  \frac{\ln(n)^{3}%
}{n^{5}}+\left(  6C^{2}-13C+\frac{15}{2}\right)  \frac{\ln(n)^{2}}{n^{5}}\\
&  +\left(  4C^{3}-13C^{2}+15C-\frac{35}{6}\right)  \frac{\ln(n)}{n^{5}%
}+\left(  C^{4}-\frac{13}{3}C^{3}+\frac{15}{2}C^{2}-\frac{35}{6}C+\frac
{61}{36}\right)  \frac{1}{n^{5}}.
\end{align*}
The expansion needn't stop here:\ $P_{5}(X)$ contributes the following:%
\begin{align*}
&  -\frac{\ln(n)^{5}}{n^{6}}-\left(  5C-\frac{77}{12}\right)  \frac{\ln
(n)^{4}}{n^{6}}-\left(  10C^{2}-\frac{77}{3}C+\frac{101}{6}\right)  \frac
{\ln(n)^{3}}{n^{6}}\\
&  -\left(  10C^{3}-\frac{77}{2}C^{2}+\frac{101}{2}C-\frac{265}{12}\right)
\frac{\ln(n)^{2}}{n^{6}}-\left(  5C^{4}-\frac{77}{3}C^{3}+\frac{101}{2}%
C^{2}-\frac{265}{6}C+\frac{515}{36}\right)  \frac{\ln(n)}{n^{6}}\\
&  -\left(  C^{5}-\frac{77}{12}C^{4}+\frac{101}{6}C^{3}-\frac{256}{12}%
C^{2}+\frac{515}{36}C-\frac{2609}{720}\right)  \frac{1}{n^{6}}%
\end{align*}
and $P_{6}(X)$ is left as an exercise. \ Under the assumption that $x_{0}%
=1/2$, the constant $C$ is estimated to be
\[
C=1.76799378613615405044...
\]
by a simple numerical method \cite{F1-sinusoid} using the preceding expansion.
\ Schoenfield \cite{Sc-sinusoid} found the series for $x_{n}$ to order
$1/n^{20}$ (by a wholly different algorithm) and calculated over 1000 digits
of $C$. \ Constants like this are highly sensitive to the initial condition
and to the choice of iterated function. \ No closed-form expression for $C$
exists in a nontrivial setting, as far as is known. \ 

\section{Radical}

For the radical function $r(x)=\frac{1}{2}\left(  -1+\sqrt{1+4x}\right)  $, we
have $\tau=1$,%
\[
\{a_{m}\}_{m=1}^{7}=\left\{  -1,2,-5,14,-42,132,-429\right\}
\]
and $\lambda=1$; consequently%
\[
\{b_{j}\}_{j=1}^{6}=\left\{  -1,2,-5,14,-42,132\right\}  ,
\]%
\[
\{a_{0j}\}_{j=1}^{6}=\left\{  1,-\frac{3}{2},\frac{10}{3},-\frac{35}{4}%
,\frac{126}{5},-77\right\}  ,
\]%
\[
\{c_{i}\}_{i=1}^{5}=\left\{  \frac{1}{2},-\frac{1}{3},\frac{13}{36}%
,-\frac{113}{240},\frac{1187}{1800}\right\}  ,
\]%
\[
T_{2}=-X-\dfrac{1}{2},
\]%
\[
T_{3}=\dfrac{1}{2}X^{2}+\dfrac{3}{2}X+\dfrac{5}{6},
\]%
\[
T_{4}=-\dfrac{1}{3}X^{3}-2X^{2}-\dfrac{19}{6}X-\dfrac{13}{9},
\]%
\[
T_{5}=\dfrac{1}{4}X^{4}+\dfrac{7}{3}X^{3}+\dfrac{27}{4}X^{2}+\dfrac{15}%
{2}X+\dfrac{2009}{720},
\]%
\[
T_{6}=-\dfrac{1}{5}X^{5}-\dfrac{31}{12}X^{4}-\dfrac{34}{3}X^{3}-\dfrac{87}%
{4}X^{2}-\dfrac{3359}{180}X-\dfrac{6973}{1200}%
\]
and
\[
P_{2}=X^{2}-X+\dfrac{1}{2},
\]%
\[
P_{3}=X^{3}-\dfrac{5}{2}X^{2}+\dfrac{5}{2}X-\dfrac{5}{6},
\]%
\[
P_{4}=X^{4}-\dfrac{13}{3}X^{3}+\dfrac{15}{2}X^{2}-\dfrac{35}{6}X+\dfrac
{61}{36},
\]%
\[
P_{5}=X^{5}-\dfrac{77}{12}X^{4}+\dfrac{101}{6}X^{3}-\dfrac{265}{12}%
X^{2}+\dfrac{515}{36}X-\dfrac{2609}{720},
\]%
\[
P_{6}=X^{6}-\dfrac{87}{10}X^{5}+\dfrac{95}{3}X^{4}-61X^{3}+65X^{2}%
-\dfrac{12977}{360}X+\dfrac{29069}{3600}.
\]
In formula (A) connecting $P_{m}$ and asymptotics of $x_{n}=r(x_{n-1})$,
replace $C$ by $-C$:%
\begin{equation}
x_{n}\sim\left(  \frac{\lambda}{n}\right)  ^{1/\tau}\left\{  1+%
{\displaystyle\sum\limits_{m=1}^{6}}
P_{m}\left(  -\frac{1}{\tau}\left[  b_{1}\ln(n)-C\right]  \right)  \frac
{1}{n^{m}}\right\}  \tag{B}\label{B}%
\end{equation}
which implies
\begin{align*}
x_{n}  &  \sim\frac{1}{n}+\frac{\ln(n)}{n^{2}}+\frac{C}{n^{2}}+\frac
{\ln(n)^{2}}{n^{3}}+(2C-1)\frac{\ln(n)}{n^{3}}+\left(  C^{2}-C+\frac{1}%
{2}\right)  \frac{1}{n^{3}}\\
&  +\frac{\ln(n)^{3}}{n^{4}}+\left(  3C-\frac{5}{2}\right)  \frac{\ln(n)^{2}%
}{n^{4}}+\left(  3C^{2}-5C+\frac{5}{2}\right)  \frac{\ln(n)}{n^{4}}+\left(
C^{3}-\frac{5}{2}C^{2}+\frac{5}{2}C-\frac{5}{6}\right)  \frac{1}{n^{4}}\\
&  +\frac{\ln(n)^{4}}{n^{5}}+\left(  4C-\frac{13}{3}\right)  \frac{\ln(n)^{3}%
}{n^{5}}+\left(  6C^{2}-13C+\frac{15}{2}\right)  \frac{\ln(n)^{2}}{n^{5}}\\
&  +\left(  4C^{3}-13C^{2}+15C-\frac{35}{6}\right)  \frac{\ln(n)}{n^{5}%
}+\left(  C^{4}-\frac{13}{3}C^{3}+\frac{15}{2}C^{2}-\frac{35}{6}C+\frac
{61}{36}\right)  \frac{1}{n^{5}}.
\end{align*}
The expansion needn't stop here:\ $P_{5}(X)$ contributes the following:%
\begin{align*}
&  \frac{\ln(n)^{5}}{n^{6}}+\left(  5C-\frac{77}{12}\right)  \frac{\ln(n)^{4}%
}{n^{6}}+\left(  10C^{2}-\frac{77}{3}C+\frac{101}{6}\right)  \frac{\ln(n)^{3}%
}{n^{6}}\\
&  +\left(  10C^{3}-\frac{77}{2}C^{2}+\frac{101}{2}C-\frac{265}{12}\right)
\frac{\ln(n)^{2}}{n^{6}}+\left(  5C^{4}-\frac{77}{3}C^{3}+\frac{101}{2}%
C^{2}-\frac{265}{6}C+\frac{515}{36}\right)  \frac{\ln(n)}{n^{6}}\\
&  +\left(  C^{5}-\frac{77}{12}C^{4}+\frac{101}{6}C^{3}-\frac{256}{12}%
C^{2}+\frac{515}{36}C-\frac{2609}{720}\right)  \frac{1}{n^{6}}%
\end{align*}
and again $P_{6}(X)$ is left as an exercise.

Comparing the results for $r(x)$ against analogous results for the logistic
function $\ell(x)$, we see that $c_{n}$ are all positive for $\ell$ while
$c_{n}$ are alternating for $r$; coefficients of $T_{n}$ are alternating for
$\ell$ while coefficients of $T_{n}$ possess the same sign for $r$;
coefficients of $P_{n}$ are positive for $\ell$ while coefficients of $P_{n}$
are alternating for $r$; disregarding signs, all corresponding coefficients
are equal. \ 

Most astonishing is the following:\ when considering $\ell$, the asymptotic
expansion for $x_{n}$ has alternating blocks of positive signs and of negative
signs; while when considering $r$, the signs of all 21 terms up to order
$1/n^{6}$ are positive; disregarding signs, all corresponding coefficients are
equal. \ Define the functions $\ell$ and $r$ to be \textbf{kindred}. \ We have
seen phenomena like this before \cite{F2-sinusoid, F3-sinusoid}, although the
earlier two pairs of kindred functions involved $x_{n}\rightarrow\infty$, not
$x_{n}\rightarrow0$.

The commonality in structure between $\ell$ and $r$ does not assist in the
numerical calculation of the constant $C$. \ Under the assumption that
$x_{0}=1/2$, we have
\[
C=-2.88756384412875082823...
\]
for $r$, using the preceding expansion. \ No relationship between constants
for $\ell$ and $r$ is observed (nor was any expected). \ We shall discuss
kindred functions in greater depth later, after examining more examples.

\section{Logarithm}

For the function $f(x)=\ln\left(  1+x\right)  $, we have $\tau=1$,%
\[
\{a_{m}\}_{m=1}^{7}=\left\{  -\frac{1}{2},\frac{1}{3},-\frac{1}{4},\frac{1}%
{5},-\frac{1}{6},\frac{1}{7},-\frac{1}{8}\right\}
\]
and $\lambda=2$; consequently
\[
\{b_{j}\}_{j=1}^{6}=\left\{  -\frac{1}{3},\frac{1}{3},-\frac{19}{45},\frac
{3}{5},-\frac{863}{945},\frac{275}{189}\right\}  ,
\]%
\[
\{a_{0j}\}_{j=1}^{6}=\left\{  1,-\frac{5}{6},1,-\frac{251}{180},\frac{19}%
{9},-\frac{19087}{5670}\right\}  ,
\]%
\[
\{c_{i}\}_{i=1}^{5}=\left\{  \frac{1}{18},-\frac{1}{135},-\frac{1}{972}%
,\frac{71}{27216},-\frac{8759}{5103000}\right\}  ,
\]%
\[
T_{2}=-\frac{1}{3}X-\dfrac{1}{18},
\]%
\[
T_{3}=\dfrac{1}{6}X^{2}+\dfrac{1}{6}X+\dfrac{7}{270},
\]%
\[
T_{4}=-\dfrac{1}{9}X^{3}-\frac{2}{9}X^{2}-\dfrac{29}{270}X-\dfrac{13}{1215}%
\]
(correcting $b_{3}$, $a_{0,4}$, $c_{2}$ values \&\ $T_{3}$ final coefficient
in \cite{ML-sinusoid})\ and
\[
P_{2}=X^{2}-\frac{1}{3}X+\dfrac{1}{18},
\]%
\[
P_{3}=X^{3}-\dfrac{5}{6}X^{2}+\dfrac{5}{18}X-\dfrac{7}{270},
\]%
\[
P_{4}=X^{4}-\dfrac{13}{9}X^{3}+\dfrac{5}{6}X^{2}-\dfrac{53}{270}X+\dfrac
{67}{4860}.
\]
In formula (A) connecting $P_{m}$ and asymptotics of $x_{n}=f(x_{n-1})$,
replace $C$ by $2C$ (for consistency's sake with \cite{Po-sinusoid}):%
\[
x_{n}\sim\left(  \frac{\lambda}{n}\right)  ^{1/\tau}\left\{  1+%
{\displaystyle\sum\limits_{m=1}^{4}}
P_{m}\left(  -\frac{1}{\tau}\left[  b_{1}\ln(n)+2C\right]  \right)  \frac
{1}{n^{m}}\right\}
\]
which implies%
\begin{align*}
x_{n}  &  \sim\frac{2}{n}+\frac{2}{3}\frac{\ln(n)}{n^{2}}-\frac{4C}{n^{2}%
}+\frac{2}{9}\frac{\ln(n)^{2}}{n^{3}}-\left(  \frac{8}{3}C+\frac{2}{9}\right)
\frac{\ln(n)}{n^{3}}+\left(  8C^{2}+\frac{4}{3}C+\frac{1}{9}\right)  \frac
{1}{n^{3}}\\
&  +\frac{2}{27}\frac{\ln(n)^{3}}{n^{4}}-\left(  \frac{4}{3}C+\frac{5}%
{27}\right)  \frac{\ln(n)^{2}}{n^{4}}+\left(  8C^{2}+\frac{20}{9}C+\frac
{5}{27}\right)  \frac{\ln(n)}{n^{4}}\\
&  -\left(  16C^{3}+\frac{20}{3}C^{2}+\frac{10}{9}C+\frac{7}{135}\right)
\frac{1}{n^{4}}+\frac{2}{81}\frac{\ln(n)^{4}}{n^{5}}-\left(  \frac{16}%
{27}C+\frac{26}{243}\right)  \frac{\ln(n)^{3}}{n^{5}}\\
&  +\left(  \frac{16}{3}C^{2}+\frac{52}{27}C+\frac{5}{27}\right)  \frac
{\ln(n)^{2}}{n^{5}}-\left(  \frac{64}{3}C^{3}+\frac{104}{9}C^{2}+\frac{20}%
{9}C+\frac{53}{405}\right)  \frac{\ln(n)}{n^{5}}\\
&  +\left(  32C^{4}+\frac{208}{9}C^{3}+\frac{20}{3}C^{2}+\frac{106}%
{135}C+\frac{67}{2430}\right)  \frac{1}{n^{5}}.
\end{align*}
Under the assumption that $x_{0}=1/2$, the constant $C$ is estimated to be%
\[
C=2.23775826599229897691...
\]
by a simple numerical method \cite{F1-sinusoid} using the preceding expansion.

\section{Exponential}

For the function $f(x)=1-\exp(-x)$, we have $\tau=1$,
\[
\{a_{m}\}_{m=1}^{7}=\left\{  -\frac{1}{2},\frac{1}{6},-\frac{1}{24},\frac
{1}{120},-\frac{1}{720},\frac{1}{5040},-\frac{1}{40320}\right\}
\]
and $\lambda=2$; consequently
\[
\{b_{j}\}_{j=1}^{6}=\left\{  \frac{1}{3},0,-\frac{1}{45},0,\frac{2}%
{945},0\right\}  ,
\]%
\[
\{a_{0j}\}_{j=1}^{6}=\left\{  1,-\frac{1}{6},0,\frac{1}{180},0,-\frac{1}%
{2835}\right\}  ,
\]%
\[
\{c_{i}\}_{i=1}^{5}=\left\{  \frac{1}{18},\frac{1}{135},-\frac{1}{972}%
,-\frac{71}{27216},-\frac{8759}{5103000}\right\}  ,
\]%
\[
T_{2}=\frac{1}{3}X-\dfrac{1}{18},
\]%
\[
T_{3}=-\dfrac{1}{6}X^{2}+\dfrac{1}{6}X-\dfrac{7}{270},
\]%
\[
T_{4}=\dfrac{1}{9}X^{3}-\frac{2}{9}X^{2}+\dfrac{29}{270}X-\dfrac{13}{1215}%
\]
and
\[
P_{2}=X^{2}+\frac{1}{3}X+\dfrac{1}{18},
\]%
\[
P_{3}=X^{3}+\frac{5}{6}X^{2}+\dfrac{5}{18}X+\dfrac{7}{270},
\]%
\[
P_{4}=X^{4}+\dfrac{13}{9}X^{3}+\dfrac{5}{6}X^{2}+\dfrac{53}{270}X+\dfrac
{67}{4860}.
\]
In formula (B) connecting $P_{m}$ and asymptotics of $x_{n}=f(x_{n-1})$,
replace $C$ by $2C$:%
\[
x_{n}\sim\left(  \frac{\lambda}{n}\right)  ^{1/\tau}\left\{  1+%
{\displaystyle\sum\limits_{m=1}^{4}}
P_{m}\left(  -\frac{1}{\tau}\left[  b_{1}\ln(n)-2C\right]  \right)  \frac
{1}{n^{m}}\right\}
\]
which implies
\begin{align*}
x_{n}  &  \sim\frac{2}{n}-\frac{2}{3}\frac{\ln(n)}{n^{2}}+\frac{4C}{n^{2}%
}+\frac{2}{9}\frac{\ln(n)^{2}}{n^{3}}-\left(  \frac{8}{3}C+\frac{2}{9}\right)
\frac{\ln(n)}{n^{3}}+\left(  8C^{2}+\frac{4}{3}C+\frac{1}{9}\right)  \frac
{1}{n^{3}}\\
&  -\frac{2}{27}\frac{\ln(n)^{3}}{n^{4}}+\left(  \frac{4}{3}C+\frac{5}%
{27}\right)  \frac{\ln(n)^{2}}{n^{4}}-\left(  8C^{2}+\frac{20}{9}C+\frac
{5}{27}\right)  \frac{\ln(n)}{n^{4}}\\
&  +\left(  16C^{3}+\frac{20}{3}C^{2}+\frac{10}{9}C+\frac{7}{135}\right)
\frac{1}{n^{4}}+\frac{2}{81}\frac{\ln(n)^{4}}{n^{5}}-\left(  \frac{16}%
{27}C+\frac{26}{243}\right)  \frac{\ln(n)^{3}}{n^{5}}\\
&  +\left(  \frac{16}{3}C^{2}+\frac{52}{27}C+\frac{5}{27}\right)  \frac
{\ln(n)^{2}}{n^{5}}-\left(  \frac{64}{3}C^{3}+\frac{104}{9}C^{2}+\frac{20}%
{9}C+\frac{53}{405}\right)  \frac{\ln(n)}{n^{5}}\\
&  +\left(  32C^{4}+\frac{208}{9}C^{3}+\frac{20}{3}C^{2}+\frac{106}%
{135}C+\frac{67}{2430}\right)  \frac{1}{n^{5}}.
\end{align*}
Observations in Section 4 carry over here: $\ln\left(  1+x\right)  $ and
$1-\exp(-x)$ are kindred functions. \ Under the assumption that $x_{0}=1/2$,
the constant $C$ is estimated to be
\[
C=-1.77611295395085782901...
\]
by a simple numerical method \cite{F1-sinusoid} using the preceding expansion.

\section{Circular Sine}

For the function $f(x)=\sin(x)$, we have $\tau=2$,
\[
\{a_{m}\}_{m=1}^{7}=\left\{  -\frac{1}{6},\frac{1}{120},-\frac{1}{5040}%
,\frac{1}{362880},-\frac{1}{39916800},\frac{1}{6227020800},-\frac
{1}{1307674368000}\right\}
\]
and $\lambda=3$; consequently
\[
\{b_{j}\}_{j=1}^{6}=\left\{  \frac{3}{5},\frac{2}{7},\frac{3}{25},\frac
{18}{385},\frac{1382}{79625},\frac{12}{1925}\right\}  ,
\]%
\[
\{a_{0j}\}_{j=1}^{6}=\left\{  1,\frac{1}{10},\frac{2}{105},\frac{3}{700}%
,\frac{2}{1925},\frac{691}{2627625}\right\}  ,
\]%
\[
\{c_{i}\}_{i=1}^{5}=\left\{  \frac{79}{350},\frac{87}{875},\frac
{91543}{1347500},\frac{18222899}{350350000},\frac{88627739}{2358125000}%
\right\}  ,
\]%
\[
T_{2}=\frac{3}{5}X-\dfrac{79}{350},
\]%
\[
T_{3}=-\dfrac{3}{10}X^{2}+\dfrac{41}{70}X-\dfrac{411}{1750},
\]%
\[
T_{4}=\dfrac{1}{5}X^{3}-\frac{134}{175}X^{2}+\dfrac{1437}{1750}X-\dfrac
{87519}{336875}%
\]
and
\[
P_{2}=\frac{3}{2}X^{2}+\frac{3}{5}X+\dfrac{79}{700},
\]%
\[
P_{3}=\frac{5}{2}X^{3}+\frac{12}{5}X^{2}+\dfrac{647}{700}X+\dfrac{411}{3500},
\]%
\[
P_{4}=\frac{35}{8}X^{4}+\dfrac{71}{10}X^{3}+\dfrac{187}{40}X^{2}+\dfrac
{2409}{1750}X+\dfrac{1606257}{10780000}.
\]
From formula (A) connecting $P_{m}$ and asymptotics of $x_{n}=f(x_{n-1})$, we
deduce
\begin{align*}
\frac{x_{n}}{\sqrt{3}}  &  \sim\frac{1}{n^{1/2}}-\frac{3}{10}\frac{\ln
(n)}{n^{3/2}}-\frac{C}{2}\frac{1}{n^{3/2}}+\frac{27}{200}\frac{\ln(n)^{2}%
}{n^{5/2}}+\left(  \frac{9}{20}C-\frac{9}{50}\right)  \frac{\ln(n)}{n^{5/2}}\\
&  +\left(  \frac{3}{8}C^{2}-\frac{3}{10}C+\frac{79}{700}\right)  \frac
{1}{n^{5/2}}-\frac{27}{400}\frac{\ln(n)^{3}}{n^{7/2}}-\left(  \frac{27}%
{80}C-\frac{27}{125}\right)  \frac{\ln(n)^{2}}{n^{7/2}}\\
&  -\left(  \frac{9}{16}C^{2}-\frac{18}{25}C+\frac{1941}{7000}\right)
\frac{\ln(n)}{n^{7/2}}-\left(  \frac{5}{16}C^{3}-\frac{3}{5}C^{2}+\frac
{647}{1400}C-\frac{411}{3500}\right)  \frac{1}{n^{7/2}}\\
&  +\frac{567}{16000}\frac{\ln(n)^{4}}{n^{9/2}}+\left(  \frac{189}{800}%
C-\frac{1917}{10000}\right)  \frac{\ln(n)^{3}}{n^{9/2}}+\left(  \frac
{189}{320}C^{2}-\frac{1917}{2000}C+\frac{1683}{4000}\right)  \frac{\ln(n)^{2}%
}{n^{9/2}}\\
&  +\left(  \frac{21}{32}C^{3}-\frac{639}{400}C^{2}+\frac{561}{400}%
C-\frac{7227}{17500}\right)  \frac{\ln(n)}{n^{9/2}}\\
&  +\left(  \frac{35}{128}C^{4}-\frac{71}{80}C^{3}+\frac{187}{160}C^{2}%
-\frac{2409}{3500}C+\frac{1606257}{10780000}\right)  \frac{1}{n^{9/2}}.
\end{align*}
This extends results provided by \cite{dB-sinusoid, BR-sinusoid} and corrects
a series given in \cite{IS-sinusoid}. \ Under the assumption that $x_{0}%
=\pi/2$, the constant $C$ is estimated to be%
\[
C=1.43045534652867724470...
\]
and $\pi/2$ is the argument at which $C$ is minimal. \ At arguments $\pi
/3,\pi/4$, $\pi/6$ we have values%
\[%
\begin{array}
[c]{ccccc}%
2.23217214236864952692..., &  & 3.96568516776811188899..., &  &
9.52859799064212800035...
\end{array}
\]
respectively. \ A popular version of this problem \cite{Sp-sinusoid} involves
the initial condition $y_{0}=1=x_{1}$ and, more generally, $y_{n}=x_{n+1}$.
\ We have%
\[
\frac{y_{n}}{\sqrt{3}}\sim\frac{1}{n^{1/2}}-\frac{3}{10}\frac{\ln(n)}{n^{3/2}%
}-\frac{C_{y}}{2}\frac{1}{n^{3/2}},
\]%
\begin{align*}
\frac{x_{n+1}}{\sqrt{3}}  &  \sim\frac{1}{(n+1)^{1/2}}-\frac{3}{10}\frac
{\ln(n+1)}{(n+1)^{3/2}}-\frac{C_{x}}{2}\frac{1}{(n+1)^{3/2}}\\
&  \sim\frac{1}{n^{1/2}}-\frac{1}{2}\,\frac{1-\frac{3}{4n}}{n^{3/2}}-\frac
{3}{10}\frac{\left(  1-\frac{3}{2n}\right)  \ln(n)}{n^{3/2}}-\frac{C_{x}}%
{2}\,\frac{1-\frac{3}{2n}}{n^{3/2}}%
\end{align*}
thus $0=-C_{y}+1+C_{x}$ upon taking differences, and $C_{y}%
=2.43045534652867724470...$.

\section{Inverse Hyperbolic Sine}

For the function $f(x)=\operatorname{arcsinh}(x)=\ln\left(  x+\sqrt{1+x^{2}%
}\right)  $, we have $\tau=2$,
\[
\{a_{m}\}_{m=1}^{7}=\left\{  -\frac{1}{6},\frac{3}{40},-\frac{5}{112}%
,\frac{35}{1152},-\frac{63}{2816},\frac{231}{13312},-\frac{143}{10240}%
\right\}
\]
and $\lambda=3$; consequently
\[
\{b_{j}\}_{j=1}^{6}=\left\{  -\frac{3}{5},\frac{31}{35},-\frac{289}{175}%
,\frac{951}{275},-\frac{6803477}{875875},\frac{3203699}{175175}\right\}  ,
\]%
\[
\{a_{0j}\}_{j=1}^{6}=\left\{  1,-\frac{11}{10},\frac{191}{105},-\frac
{2497}{700},\frac{14797}{1925},-\frac{92427157}{5255250}\right\}  ,
\]%
\[
\{c_{i}\}_{i=1}^{5}=\left\{  \frac{79}{350},-\frac{87}{875},\frac
{91543}{1347500},-\frac{18222899}{350350000},\frac{88627739}{2358125000}%
\right\}  ,
\]%
\[
T_{2}=-\frac{3}{5}X-\dfrac{79}{350},
\]%
\[
T_{3}=\dfrac{3}{10}X^{2}+\dfrac{41}{70}X+\dfrac{411}{1750},
\]%
\[
T_{4}=-\dfrac{1}{5}X^{3}-\frac{134}{175}X^{2}-\dfrac{1437}{1750}%
X-\dfrac{87519}{336875}%
\]
(correcting $b_{3}$, $b_{4}$ $a_{0,4}$, $c_{2}$ values \&\ $T_{3}$ final
coefficient in \cite{ML-sinusoid})\ and%
\[
P_{2}=\frac{3}{2}X^{2}-\frac{3}{5}X+\dfrac{79}{700},
\]%
\[
P_{3}=\frac{5}{2}X^{3}-\frac{12}{5}X^{2}+\dfrac{647}{700}X-\dfrac{411}{3500},
\]%
\[
P_{4}=\frac{35}{8}X^{4}-\dfrac{71}{10}X^{3}+\dfrac{187}{40}X^{2}-\dfrac
{2409}{1750}X+\dfrac{1606257}{10780000}.
\]
From formula (B) connecting $P_{m}$ and asymptotics of $x_{n}=f(x_{n-1})$, we deduce%

\begin{align*}
\frac{x_{n}}{\sqrt{3}}  &  \sim\frac{1}{n^{1/2}}+\frac{3}{10}\frac{\ln
(n)}{n^{3/2}}+\frac{C}{2}\frac{1}{n^{3/2}}+\frac{27}{200}\frac{\ln(n)^{2}%
}{n^{5/2}}+\left(  \frac{9}{20}C-\frac{9}{50}\right)  \frac{\ln(n)}{n^{5/2}}\\
&  +\left(  \frac{3}{8}C^{2}-\frac{3}{10}C+\frac{79}{700}\right)  \frac
{1}{n^{5/2}}+\frac{27}{400}\frac{\ln(n)^{3}}{n^{7/2}}+\left(  \frac{27}%
{80}C-\frac{27}{125}\right)  \frac{\ln(n)^{2}}{n^{7/2}}\\
&  +\left(  \frac{9}{16}C^{2}-\frac{18}{25}C+\frac{1941}{7000}\right)
\frac{\ln(n)}{n^{7/2}}+\left(  \frac{5}{16}C^{3}-\frac{3}{5}C^{2}+\frac
{647}{1400}C-\frac{411}{3500}\right)  \frac{1}{n^{7/2}}\\
&  +\frac{567}{16000}\frac{\ln(n)^{4}}{n^{9/2}}+\left(  \frac{189}{800}%
C-\frac{1917}{10000}\right)  \frac{\ln(n)^{3}}{n^{9/2}}+\left(  \frac
{189}{320}C^{2}-\frac{1917}{2000}C+\frac{1683}{4000}\right)  \frac{\ln(n)^{2}%
}{n^{9/2}}\\
&  +\left(  \frac{21}{32}C^{3}-\frac{639}{400}C^{2}+\frac{561}{400}%
C-\frac{7227}{17500}\right)  \frac{\ln(n)}{n^{9/2}}\\
&  +\left(  \frac{35}{128}C^{4}-\frac{71}{80}C^{3}+\frac{187}{160}C^{2}%
-\frac{2409}{3500}C+\frac{1606257}{10780000}\right)  \frac{1}{n^{9/2}}.
\end{align*}
Observations in Section 4 carry over here: $\sin\left(  x\right)  $ and
$\operatorname{arcsinh}(x)$ are kindred functions (answering the question in
this paper's title). \ For reasons of space, we stop here.

\section{Inverse Circular Tangent}

For the function $f(x)=\arctan(x)$, we have $\tau=2$,
\[
\{a_{m}\}_{m=1}^{7}=\left\{  -\frac{1}{3},\frac{1}{5},-\frac{1}{7},\frac{1}%
{9},-\frac{1}{11},\frac{1}{13},-\frac{1}{15}\right\}
\]
and $\lambda=3/2$; consequently
\[
\{b_{j}\}_{j=1}^{6}=\left\{  -\frac{3}{20},\frac{4}{35},-\frac{19}{175}%
,\frac{222}{1925},-\frac{459257}{3503500},\frac{109271}{700700}\right\}  ,
\]%
\[
\{a_{0j}\}_{j=1}^{6}=\left\{  1,-\frac{13}{20},\frac{251}{420},-\frac
{3551}{5600},\frac{22417}{30800},-\frac{147636491}{168168000}\right\}  ,
\]%
\[
\{c_{i}\}_{i=1}^{5}=\left\{  \frac{47}{2800},\frac{3}{16000},-\frac
{83723}{86240000},\frac{27832729}{89689600000},\frac{2800730629}%
{15695680000000}\right\}  ,
\]

\[
T_{2}=-\frac{3}{20}X-\frac{47}{2800},
\]%
\[
T_{3}=\frac{3}{40}X^{2}+\frac{11}{280}X+\frac{261}{112000},
\]%
\[
T_{4}=-\frac{1}{20}X^{3}-\frac{283}{5600}X^{2}-\frac{591}{56000}X+\frac
{58557}{172480000}%
\]
and
\[
P_{2}=\frac{3}{2}X^{2}-\frac{3}{20}X+\dfrac{47}{5600},
\]%
\[
P_{3}=\frac{5}{2}X^{3}-\frac{3}{5}X^{2}+\dfrac{361}{5600}X-\dfrac{261}%
{224000},
\]%
\[
P_{4}=\frac{35}{8}X^{4}-\dfrac{71}{40}X^{3}+\dfrac{101}{320}X^{2}-\dfrac
{3993}{224000}X-\dfrac{44217}{689920000}.
\]
From formula (B) connecting $P_{m}$ and asymptotics of $x_{n}=f(x_{n-1})$, we deduce%

\begin{align*}
\frac{x_{n}}{\sqrt{3/2}}  &  \sim\frac{1}{n^{1/2}}+\frac{3}{40}\frac{\ln
(n)}{n^{3/2}}+\frac{C}{2}\frac{1}{n^{3/2}}+\frac{27}{3200}\frac{\ln(n)^{2}%
}{n^{5/2}}+\left(  \frac{9}{80}C-\frac{9}{800}\right)  \frac{\ln(n)}{n^{5/2}%
}\\
&  +\left(  \frac{3}{8}C^{2}-\frac{3}{40}C+\frac{47}{5600}\right)  \frac
{1}{n^{5/2}}+\frac{27}{25600}\frac{\ln(n)^{3}}{n^{7/2}}+\left(  \frac
{27}{1280}C-\frac{27}{8000}\right)  \frac{\ln(n)^{2}}{n^{7/2}}\\
&  +\left(  \frac{9}{64}C^{2}-\frac{9}{200}C+\frac{1083}{224000}\right)
\frac{\ln(n)}{n^{7/2}}+\left(  \frac{5}{16}C^{3}-\frac{3}{20}C^{2}+\frac
{361}{11200}C-\frac{261}{224000}\right)  \frac{1}{n^{7/2}}\\
&  +\frac{567}{4096000}\frac{\ln(n)^{4}}{n^{9/2}}+\left(  \frac{189}%
{51200}C-\frac{1917}{2560000}\right)  \frac{\ln(n)^{3}}{n^{9/2}}+\left(
\frac{189}{5120}C^{2}-\frac{1917}{128000}C+\frac{909}{512000}\right)
\frac{\ln(n)^{2}}{n^{9/2}}\\
&  +\left(  \frac{21}{128}C^{3}-\frac{639}{6400}C^{2}+\frac{303}{12800}%
C-\frac{11979}{8960000}\right)  \frac{\ln(n)}{n^{9/2}}\\
&  +\left(  \frac{35}{128}C^{4}-\frac{71}{320}C^{3}+\frac{101}{1280}%
C^{2}-\frac{3993}{448000}C-\frac{44217}{689920000}\right)  \frac{1}{n^{9/2}}.
\end{align*}

\section{Hyperbolic Tangent}

For the function $f(x)=\tanh(x)=[\exp(2x)-1]/[\exp(2x)+1]$, we have $\tau=2$,
\[
\{a_{m}\}_{m=1}^{7}=\left\{  -\frac{1}{3},\frac{2}{15},-\frac{17}{315}%
,\frac{62}{2835},-\frac{1382}{155925},\frac{21844}{6081075},-\frac
{929569}{638512875}\right\}
\]
and $\lambda=3/2$; consequently
\[
\{b_{j}\}_{j=1}^{6}=\left\{  \frac{3}{20},-\frac{1}{28},\frac{3}{400}%
,-\frac{9}{6160},\frac{691}{2548000},-\frac{3}{61600}\right\}  ,
\]%
\[
\{a_{0j}\}_{j=1}^{6}=\left\{  1,-\frac{7}{20},\frac{31}{210},-\frac{381}%
{5600},\frac{73}{2200},-\frac{1414477}{84084000}\right\}  ,
\]%
\[
\{c_{i}\}_{i=1}^{5}=\left\{  \frac{47}{2800},-\frac{3}{16000},-\frac
{83723}{86240000},-\frac{27832729}{89689600000},\frac{2800730629}%
{15695680000000}\right\}  ,
\]%
\[
T_{2}=\frac{3}{20}X-\frac{47}{2800},
\]%
\[
T_{3}=-\frac{3}{40}X^{2}+\frac{11}{280}X-\frac{261}{112000},
\]%
\[
T_{4}=\frac{1}{20}X^{3}-\frac{283}{5600}X^{2}+\frac{591}{56000}X+\frac
{58557}{172480000}%
\]
($c_{2}$ in \cite{ML-sinusoid} should be multipled by $1/4$)\ and%
\[
P_{2}=\frac{3}{2}X^{2}+\frac{3}{20}X+\dfrac{47}{5600},
\]%
\[
P_{3}=\frac{5}{2}X^{3}+\frac{3}{5}X^{2}+\dfrac{361}{5600}X+\dfrac{261}%
{224000},
\]%
\[
P_{4}=\frac{35}{8}X^{4}+\dfrac{71}{40}X^{3}+\dfrac{101}{320}X^{2}+\dfrac
{3993}{224000}X-\dfrac{44217}{689920000}.
\]
From formula (A) connecting $P_{m}$ and asymptotics of $x_{n}=f(x_{n-1})$, we deduce%

\begin{align*}
\frac{x_{n}}{\sqrt{3/2}}  &  \sim\frac{1}{n^{1/2}}-\frac{3}{40}\frac{\ln
(n)}{n^{3/2}}-\frac{C}{2}\frac{1}{n^{3/2}}+\frac{27}{3200}\frac{\ln(n)^{2}%
}{n^{5/2}}+\left(  \frac{9}{80}C-\frac{9}{800}\right)  \frac{\ln(n)}{n^{5/2}%
}\\
&  +\left(  \frac{3}{8}C^{2}-\frac{3}{40}C+\frac{47}{5600}\right)  \frac
{1}{n^{5/2}}-\frac{27}{25600}\frac{\ln(n)^{3}}{n^{7/2}}-\left(  \frac
{27}{1280}C-\frac{27}{8000}\right)  \frac{\ln(n)^{2}}{n^{7/2}}\\
&  -\left(  \frac{9}{64}C^{2}-\frac{9}{200}C+\frac{1083}{224000}\right)
\frac{\ln(n)}{n^{7/2}}-\left(  \frac{5}{16}C^{3}-\frac{3}{20}C^{2}+\frac
{361}{11200}C-\frac{261}{224000}\right)  \frac{1}{n^{7/2}}\\
&  +\frac{567}{4096000}\frac{\ln(n)^{4}}{n^{9/2}}+\left(  \frac{189}%
{51200}C-\frac{1917}{2560000}\right)  \frac{\ln(n)^{3}}{n^{9/2}}+\left(
\frac{189}{5120}C^{2}-\frac{1917}{128000}C+\frac{909}{512000}\right)
\frac{\ln(n)^{2}}{n^{9/2}}\\
&  +\left(  \frac{21}{128}C^{3}-\frac{639}{6400}C^{2}+\frac{303}{12800}%
C-\frac{11979}{8960000}\right)  \frac{\ln(n)}{n^{9/2}}\\
&  +\left(  \frac{35}{128}C^{4}-\frac{71}{320}C^{3}+\frac{101}{1280}%
C^{2}-\frac{3993}{448000}C-\frac{44217}{689920000}\right)  \frac{1}{n^{9/2}}.
\end{align*}
Observations in Section 4 carry over here: $\arctan\left(  x\right)  $ and
$\tanh(x)$ are kindred functions.

\section{Fresnel Cosine}

For the function%
\[
f(x)=%
{\displaystyle\int\limits_{0}^{x}}
\cos\left(  \frac{\pi}{2}\,t^{2}\right)  dt
\]
we have $\tau=4$,
\[
\{a_{m}\}_{m=1}^{6}=\left\{  -\frac{\pi^{2}}{40},\frac{\pi^{4}}{3456}%
,-\frac{\pi^{6}}{599040},\frac{\pi^{8}}{175472640},-\frac{\pi^{10}%
}{78033715200},\frac{\pi^{12}}{49049763840000}\right\}
\]
and $\lambda=10/\pi^{2}$; consequently
\[
\{b_{j}\}_{j=1}^{6}=\left\{  \frac{55}{108},\frac{245}{1404},\frac
{403315}{9022104},\frac{188719}{21051576},\frac{2770688381}{1950723238464}%
,\frac{1382689765}{8081567702208}\right\}  ,
\]%
\[
\{a_{0j}\}_{j=1}^{6}=\left\{  1,\frac{1}{108},-\frac{1}{702},-\frac
{7639}{36088416},-\frac{5587}{315773640},-\frac{2445241}{2926084857696}%
\right\}  ,
\]%
\[
\{c_{i}\}_{i=1}^{5}=\left\{  \frac{25745}{151632},\frac{250903475}%
{3897548928},\frac{124092693449}{3379862723328},\frac{2012588028078415}%
{91012943413776384},\frac{168154016370424766095}{15749152767150106816512}%
\right\}  ,
\]%
\[
T_{2}=\frac{55}{108}X-\frac{25745}{151632},
\]%
\[
T_{3}=-\frac{55}{216}X^{2}+\frac{1205}{2808}X-\frac{587905525}{3897548928},
\]%
\[
T_{4}=\frac{55}{324}X^{3}-\frac{169465}{303264}X^{2}+\frac{1013788675}%
{1948774464}X-\frac{16359252550091}{114915332593152}%
\]
(correcting $b_{3}$, $a_{0,4}$, $c_{1}$, $c_{2}$ values \&\ $T_{3}$ final two
coefficients in \cite{ML-sinusoid})\ and%

\[
P_{2}=\frac{5}{2}X^{2}+\frac{55}{108}X+\dfrac{25745}{606528},
\]%
\[
P_{3}=\frac{15}{2}X^{3}+\frac{385}{108}X^{2}+\dfrac{389005}{606528}%
X+\dfrac{587905525}{15590195712},
\]%
\[
P_{4}=\frac{195}{8}X^{4}+\dfrac{12485}{648}X^{3}+\dfrac{62045}{10368}%
X^{2}+\dfrac{12734847275}{15590195712}X+\dfrac{73718758987739}%
{1838645321490432}.
\]
Let $\kappa=\sqrt{5/8}/\pi$. \ From formula (B) connecting $P_{m}$ and
asymptotics of $x_{n}=f(x_{n-1})$, we deduce
\begin{align*}
\frac{\,x_{n}}{\sqrt{\kappa}}  &  \sim\frac{2}{n^{1/4}}-\frac{55}{216}%
\frac{\ln(n)}{n^{5/4}}+\frac{C}{2}\frac{1}{n^{5/4}}+\frac{15125}{186624}%
\frac{\ln(n)^{2}}{n^{9/4}}-\left(  \frac{275}{864}C+\frac{3025}{23328}\right)
\frac{\ln(n)}{n^{9/4}}\\
&  +\left(  \frac{5}{16}C^{2}+\frac{55}{216}C+\frac{25745}{303264}\right)
\frac{1}{n^{9/4}}-\frac{831875}{26873856}\frac{\ln(n)^{3}}{n^{13/4}}+\left(
\frac{15125}{82944}C+\frac{1164625}{10077696}\right)  \frac{\ln(n)^{2}%
}{n^{13/4}}\\
&  -\left(  \frac{275}{768}C^{2}+\frac{21175}{46656}C+\frac{21395275}%
{131010048}\right)  \frac{\ln(n)}{n^{13/4}}+\left(  \frac{15}{64}C^{3}%
+\frac{385}{864}C^{2}+\frac{389005}{1213056}C+\frac{587905525}{7795097856}%
\right)  \frac{1}{n^{13/4}}\\
&  +\frac{594790625}{46438023168}\frac{\ln(n)^{4}}{n^{17/4}}-\left(
\frac{10814375}{107495424}C+\frac{2077191875}{26121388032}\right)  \frac
{\ln(n)^{3}}{n^{17/4}}\\
&  +\left(  \frac{196625}{663552}C^{2}+\frac{37767125}{80621568}%
C+\frac{187686125}{967458816}\right)  \frac{\ln(n)^{2}}{n^{17/4}}\\
&  -\left(  \frac{3575}{9216}C^{3}+\frac{686675}{746496}C^{2}+\frac
{3412475}{4478976}C+\frac{700416600125}{3367482273792}\right)  \frac{\ln
(n)}{n^{17/4}}\\
&  +\left(  \frac{195}{1024}C^{4}+\frac{12485}{20736}C^{3}+\frac{62045}%
{82944}C^{2}+\frac{12734847275}{31180391424}C+\frac{73718758987739}%
{919322660745216}\right)  \frac{1}{n^{17/4}}.
\end{align*}

\section{Kindredness I}

By now, the reader has detected that kindred functions resemble inverses,
tailored vaguely so that $a_{1}<0$. \ The inverse of $\ell(x)=x(1-x)$,
restricted to $[0,1/2]$, is
\[
\frac{1}{2}\left(  1-\sqrt{1-4x}\right)  =x+x^{2}+2x^{3}+5x^{4}+14x^{5}%
+42x^{6}+\cdots
\]
whereas%
\[
r(x)=\frac{1}{2}\left(  -1+\sqrt{1+4x}\right)  =x-x^{2}+2x^{3}-5x^{4}%
+14x^{5}-42x^{6}+-\cdots.
\]
The inverse of $\ln(1+x)$ is
\[
-1+\exp(x)=x+\frac{1}{2}x^{2}+\frac{1}{6}x^{3}+\frac{1}{24}x^{4}+\frac{1}%
{120}x^{5}+\frac{1}{720}x^{6}+\cdots
\]
whereas%
\[
1-\exp(-x)=x-\frac{1}{2}x^{2}+\frac{1}{6}x^{3}-\frac{1}{24}x^{4}+\frac{1}%
{120}x^{5}-\frac{1}{720}x^{6}+-\cdots.
\]
The series for $\operatorname{arcsinh}(x)$ is the same as $\arcsin(x)$, with
alternating signs. \ The series for $\tan(x)$ is the same as $\tanh(x)$, with
all positive signs. \ A pattern is now clear. \ While the inverse of the
Fresnel cosine $f(x)$ does not have a name, we may use Lagrange inversion to
revert the series, yielding:%

\[
x+\frac{\pi^{2}}{40}x^{5}+\frac{49\pi^{4}}{17280}x^{9}+\frac{4019\pi^{6}%
}{8985600}x^{13}+\frac{42037157\pi^{8}}{513257472000}x^{17}+\frac
{21129194183\pi^{10}}{1293408829440000}x^{21}+\cdots
\]
and then define
\[
g(x)=x-\frac{\pi^{2}}{40}x^{5}+\frac{49\pi^{4}}{17280}x^{9}-\frac{4019\pi^{6}%
}{8985600}x^{13}+\frac{42037157\pi^{8}}{513257472000}x^{17}-\frac
{21129194183\pi^{10}}{1293408829440000}x^{21}+-\cdots
\]
to be the iterated function. \ Continuing as before, we have $\tau=4$ and
$\lambda=10/\pi^{2}$; consequently%

\[
\{b_{j}\}_{j=1}^{6}=\left\{  -\frac{55}{108},\frac{80}{117},-\frac
{1613305}{1388016},\frac{343709}{154791},-\frac{181191932495}{39810678336}%
,\frac{7036909925}{719512794}\right\}  ,
\]%
\[
\{a_{0j}\}_{j=1}^{6}=\left\{  1,-\frac{109}{108},\frac{2143}{1404}%
,-\frac{98701909}{36088416},\frac{3497783}{649740},-\frac{65830792020265}%
{5852169715392}\right\}  ,
\]%
\[
\{c_{i}\}_{i=1}^{5}=\left\{  \frac{25745}{151632},-\frac{250903475}%
{3897548928},\frac{124092693449}{3379862723328},-\frac{2012588028078415}%
{91012943413776384},\frac{168154016370424766095}{15749152767150106816512}%
\right\}  ,
\]%
\[
T_{2}=-\frac{55}{108}X-\frac{25745}{151632},
\]%
\[
T_{3}=\frac{55}{216}X^{2}+\frac{1205}{2808}X+\frac{587905525}{3897548928},
\]%
\[
T_{4}=-\frac{55}{324}X^{3}-\frac{169465}{303264}X^{2}-\frac{1013788675}%
{1948774464}X-\frac{16359252550091}{114915332593152}%
\]
and%

\[
P_{2}=\frac{5}{2}X^{2}-\frac{55}{108}X+\dfrac{25745}{606528},
\]%
\[
P_{3}=\frac{15}{2}X^{3}-\frac{385}{108}X^{2}+\dfrac{389005}{606528}%
X-\dfrac{587905525}{15590195712},
\]%
\[
P_{4}=\frac{195}{8}X^{4}-\dfrac{12485}{648}X^{3}+\dfrac{62045}{10368}%
X^{2}-\dfrac{12734847275}{15590195712}X+\dfrac{73718758987739}%
{1838645321490432}.
\]
Let $\kappa=\sqrt{5/8}/\pi$. \ From formula (A) connecting $P_{m}$ and
asymptotics of $x_{n}=f(x_{n-1})$, we deduce
\begin{align*}
\frac{\,x_{n}}{\sqrt{\kappa}}  &  \sim\frac{2}{n^{1/4}}+\frac{55}{216}%
\frac{\ln(n)}{n^{5/4}}-\frac{C}{2}\frac{1}{n^{5/4}}+\frac{15125}{186624}%
\frac{\ln(n)^{2}}{n^{9/4}}-\left(  \frac{275}{864}C+\frac{3025}{23328}\right)
\frac{\ln(n)}{n^{9/4}}\\
&  +\left(  \frac{5}{16}C^{2}+\frac{55}{216}C+\frac{25745}{303264}\right)
\frac{1}{n^{9/4}}+\frac{831875}{26873856}\frac{\ln(n)^{3}}{n^{13/4}}-\left(
\frac{15125}{82944}C+\frac{1164625}{10077696}\right)  \frac{\ln(n)^{2}%
}{n^{13/4}}\\
&  +\left(  \frac{275}{768}C^{2}+\frac{21175}{46656}C+\frac{21395275}%
{131010048}\right)  \frac{\ln(n)}{n^{13/4}}-\left(  \frac{15}{64}C^{3}%
+\frac{385}{864}C^{2}+\frac{389005}{1213056}C+\frac{587905525}{7795097856}%
\right)  \frac{1}{n^{13/4}}\\
&  +\frac{594790625}{46438023168}\frac{\ln(n)^{4}}{n^{17/4}}-\left(
\frac{10814375}{107495424}C+\frac{2077191875}{26121388032}\right)  \frac
{\ln(n)^{3}}{n^{17/4}}\\
&  +\left(  \frac{196625}{663552}C^{2}+\frac{37767125}{80621568}%
C+\frac{187686125}{967458816}\right)  \frac{\ln(n)^{2}}{n^{17/4}}\\
&  -\left(  \frac{3575}{9216}C^{3}+\frac{686675}{746496}C^{2}+\frac
{3412475}{4478976}C+\frac{700416600125}{3367482273792}\right)  \frac{\ln
(n)}{n^{17/4}}\\
&  +\left(  \frac{195}{1024}C^{4}+\frac{12485}{20736}C^{3}+\frac{62045}%
{82944}C^{2}+\frac{12734847275}{31180391424}C+\frac{73718758987739}%
{919322660745216}\right)  \frac{1}{n^{17/4}}.
\end{align*}
Observations in Section 4 carry over here: $f\left(  x\right)  $ and $g(x)$
are kindred functions. \ Three recurrences given by cubic functions appeared
in \cite{F2-sinusoid}. \ It is possible to construct kindred functions for
each of these, either by use of Tartaglia's (quite complicated) formula or via
(relatively elegant) series reversion.

\section{Lambert $W$}

For the function $W(x)$ satisfying $W\exp(W)=x$, we have $\tau$ $=1$,
\[
\{a_{m}\}_{m=1}^{8}=\left\{  -1,\frac{3}{2},-\frac{8}{3},\frac{125}{24}%
,-\frac{54}{5},\frac{16807}{720},-\frac{16384}{315},\frac{531441}%
{4480}\right\}
\]
and $\lambda=1$; consequently
\[
\{b_{j}\}_{j=1}^{6}=\left\{  -\frac{1}{2},\frac{2}{3},-\frac{9}{8},\frac
{32}{15},-\frac{625}{144},\frac{324}{35}\right\}  ,
\]%
\[
\{a_{0j}\}_{j=1}^{6}=\left\{  1,-1,\frac{3}{2},-\frac{8}{3},\frac{125}%
{24},-\frac{54}{5}\right\}  ,
\]%
\[
\{c_{i}\}_{i=1}^{5}=\left\{  \frac{1}{6},-\frac{1}{16},\frac{19}{540}%
,-\frac{1}{48},\frac{41}{4200}\right\}  ,
\]%
\[
T_{2}=-\frac{1}{2}X-\frac{1}{6},
\]%
\[
T_{3}=\frac{1}{4}X^{2}+\frac{5}{12}X+\frac{7}{48},
\]%
\[
T_{4}=-\frac{1}{6}X^{3}-\frac{13}{24}X^{2}-\frac{1}{2}X-\frac{587}{4320}%
\]
and%
\[
P_{2}=X^{2}-\frac{1}{2}X+\dfrac{1}{6},
\]%
\[
P_{3}=X^{3}-\frac{5}{4}X^{2}+\dfrac{3}{4}X-\dfrac{7}{48},
\]%
\[
P_{4}=X^{4}-\dfrac{13}{6}X^{3}+\dfrac{17}{8}X^{2}-\dfrac{23}{24}X+\dfrac
{707}{4320}.
\]
From formula (B) connecting $P_{m}$ and asymptotics of $x_{n}=W(x_{n-1})$, we deduce%

\begin{align*}
x_{n}  &  \sim\frac{1}{n}+\frac{1}{2}\frac{\ln(n)}{n^{2}}+\frac{C}{n^{2}%
}+\frac{1}{4}\frac{\ln(n)^{2}}{n^{3}}+\left(  C-\frac{1}{4}\right)  \frac
{\ln(n)}{n^{3}}+\left(  C^{2}-\frac{1}{2}C+\frac{1}{6}\right)  \frac{1}{n^{3}%
}\\
&  +\frac{1}{8}\frac{\ln(n)^{3}}{n^{4}}+\left(  \frac{3}{4}C-\frac{5}%
{16}\right)  \frac{\ln(n)^{2}}{n^{4}}+\left(  \frac{3}{2}C^{2}-\frac{5}%
{4}C+\frac{3}{8}\right)  \frac{\ln(n)}{n^{4}}+\left(  C^{3}-\frac{5}{4}%
C^{2}+\frac{3}{4}C-\frac{7}{48}\right)  \frac{1}{n^{4}}\\
&  +\frac{1}{16}\frac{\ln(n)^{4}}{n^{5}}+\left(  \frac{1}{2}C-\frac{13}%
{48}\right)  \frac{\ln(n)^{3}}{n^{5}}+\left(  \frac{3}{2}C^{2}-\frac{13}%
{8}C+\frac{17}{32}\right)  \frac{\ln(n)^{2}}{n^{5}}\\
&  +\left(  2C^{3}-\frac{13}{4}C^{2}+\frac{17}{8}C-\frac{23}{48}\right)
\frac{\ln(n)}{n^{5}}+\left(  C^{4}-\frac{13}{6}C^{3}+\frac{17}{8}C^{2}%
-\frac{23}{24}C+\frac{707}{4320}\right)  \frac{1}{n^{5}}.
\end{align*}

\section{Kindredness II}

The functional inverse of $W(x)$ is $x\exp(x)$ and thus we examine
$Z(x)=x\exp(-x)$. \ We have $\tau$ $=1$,
\[
\{a_{m}\}_{m=1}^{8}=\left\{  -1,\frac{1}{2},-\frac{1}{6},\frac{1}{24}%
,-\frac{1}{120},\frac{1}{720},-\frac{1}{5040},\frac{1}{40320}\right\}
\]
and $\lambda=1$; consequently
\[
\{b_{j}\}_{j=1}^{6}=\left\{  \frac{1}{2},\frac{1}{6},\frac{1}{24},\frac
{1}{120},\frac{1}{720},\frac{1}{5040}\right\}  ,
\]%
\[
\{a_{0j}\}_{j=1}^{6}=\left\{  1,0,0,0,0\right\}  ,
\]%
\[
\{c_{i}\}_{i=1}^{5}=\left\{  \frac{1}{6},\frac{1}{16},\frac{19}{540},\frac
{1}{48},\frac{41}{4200}\right\}  ,
\]%
\[
T_{2}=\frac{1}{2}X-\frac{1}{6},
\]%
\[
T_{3}=-\frac{1}{4}X^{2}+\frac{5}{12}X-\frac{7}{48},
\]%
\[
T_{4}=\frac{1}{6}X^{3}-\frac{13}{24}X^{2}+\frac{1}{2}X-\frac{587}{4320}%
\]
and
\[
P_{2}=X^{2}+\frac{1}{2}X+\dfrac{1}{6},
\]%
\[
P_{3}=X^{3}+\frac{5}{4}X^{2}+\dfrac{3}{4}X+\dfrac{7}{48},
\]%
\[
P_{4}=X^{4}+\dfrac{13}{6}X^{3}+\dfrac{17}{8}X^{2}+\dfrac{23}{24}X+\dfrac
{707}{4320}.
\]
From formula (A) connecting $P_{m}$ and asymptotics of $x_{n}=Z(x_{n-1})$, we deduce%

\begin{align*}
x_{n}  &  \sim\frac{1}{n}-\frac{1}{2}\frac{\ln(n)}{n^{2}}-\frac{C}{n^{2}%
}+\frac{1}{4}\frac{\ln(n)^{2}}{n^{3}}+\left(  C-\frac{1}{4}\right)  \frac
{\ln(n)}{n^{3}}+\left(  C^{2}-\frac{1}{2}C+\frac{1}{6}\right)  \frac{1}{n^{3}%
}\\
&  -\frac{1}{8}\frac{\ln(n)^{3}}{n^{4}}-\left(  \frac{3}{4}C-\frac{5}%
{16}\right)  \frac{\ln(n)^{2}}{n^{4}}-\left(  \frac{3}{2}C^{2}-\frac{5}%
{4}C+\frac{3}{8}\right)  \frac{\ln(n)}{n^{4}}-\left(  C^{3}-\frac{5}{4}%
C^{2}+\frac{3}{4}C-\frac{7}{48}\right)  \frac{1}{n^{4}}\\
&  +\frac{1}{16}\frac{\ln(n)^{4}}{n^{5}}+\left(  \frac{1}{2}C-\frac{13}%
{48}\right)  \frac{\ln(n)^{3}}{n^{5}}+\left(  \frac{3}{2}C^{2}-\frac{13}%
{8}C+\frac{17}{32}\right)  \frac{\ln(n)^{2}}{n^{5}}\\
&  +\left(  2C^{3}-\frac{13}{4}C^{2}+\frac{17}{8}C-\frac{23}{48}\right)
\frac{\ln(n)}{n^{5}}+\left(  C^{4}-\frac{13}{6}C^{3}+\frac{17}{8}C^{2}%
-\frac{23}{24}C+\frac{707}{4320}\right)  \frac{1}{n^{5}}.
\end{align*}
Observations in Section 4 carry over here: $W(x)$ and $Z(x)$ are kindred
functions. \ Under the assumption that $x_{0}=1$, the constant $C$ is
estimated to be%
\[
C=-1.29024720868776429166...
\]
by a simple numerical method \cite{F1-sinusoid} using the preceding expansion.

All known examples of kindred functions involve pairs. \ Open questions
abound. \ Do kindred triples exist?\ \ Do kindred quadruples exist? \ While a
pattern for assembling a kindred pair is evident, a rigorous proof is not
clear. \ Why should tailoring a power series for an inverse function
(replacing positive signs by alternating signs and vice versa) be sufficient
to produce such vivid commonality in structure (in the asymptotic expansion of
$x_{n}$)? \ We repeat a query from \cite{F2-sinusoid}. Among nontrivial cases,
does there exist an explicit recurrence $f$ and an explicit value $x_{0}$ such
that the corresponding constant $C$ is known exactly (i.e., possesses a
closed-form expression)? \ The method in \cite{BR-sinusoid, ML-sinusoid}
conceivably could assist in finding an answer.

\section{Acknowledgements}

The creators of Mathematica earn my gratitude every day:\ this paper could not
have otherwise been written. \ I\ thank Vichian Laohakosol \cite{ML-sinusoid}
for his encouraging words.

\end{document}